\input amstex
\input epsf
\documentstyle{amsppt}
\magnification=1200
\NoBlackBoxes
\nologo
\tolerance 1000
\hsize 16.4 truecm
\vsize 23.0 truecm

\global\topskip=0pt
\baselineskip 18pt

\font\double=msbm10
\def\BbbR{\hbox{\double\char'122}}

\def\section#1{\centerline{\bf #1}}
\define\s{\sharp}

\def\L{\frak L}
\def\l{\frak l}

\def\P{\frak P}
\def\G{\frak G}

\def\ls{\leqslant}
\def\gs{\geqslant}
\def\one{\hbox{1\hskip -4pt 1}}

\def\codim{\operatorname{codim}}

\define \cl {\Cal L_{S^1}}
\define \cg {\Cal L_{\Gamma}}
\define \Ga {\Gamma}
\define \la {\lambda}
\define \pa {\partial}
\define \bR {\Bbb R}
\define \bC {\Bbb C}
\define \bZ {\Bbb F}
\define \bF {\Bbb F}
\define \al {\alpha}
\define \siij {\sigma_{i,j}}
\define \Si {\Sigma}
\define \si {\sigma}
\def\kap{\varkappa}
\def\sg{\kap}
\leftline{\eightit To V.I.Arnold}
\leftline{\eightit on the occasion of his 60-th birthday}

\vskip 3cm

\topmatter

\title
Connected components in the intersection of two open opposite 
Schubert cells in $SL_n(\bR)/B$
\endtitle 

\author 
B.~Shapiro$^{\ddag}$, M.~Shapiro$^*$ and A.~Vainshtein$^{\dag,}$
\footnote{The research of this author was supported by the Rashi
Foundation.\kern 4cm}
\endauthor

\headline{\hss{\rm\number\pageno}\hss}
\footline{}

\affil 
$^\ddag$ Department of Mathematics, University of Stockholm\\
S-10691, Sweden, {\tt shapiro\@matematik.su.se}\\ 
$^*$ Department of Mathematics, Royal Institute of Technology\\
S-10044, Sweden, {\tt mshapiro\@math.kth.se} \\ 
$^\dag$ Dept. of Mathematics and Computer Science, University of Haifa
\\ Mount Carmel, 31905 Haifa, Israel, {\tt alek\@mathcs11.haifa.ac.il} 
\endaffil


\abstract
In this paper we reduce the problem concerning the number of connected
components in the
intersection of two real opposite open Schubert cells in $SL_n(\bR)/B$
to a purely combinatorial question in the space of upper triangular 
matrices with $\bZ_2$-valued entries.
The crucial step of the reduction uses the parametrization of the space of real
unipotent totally positive upper triangular matrices  introduced in \cite
{L} and \cite{BFZ}.
\endabstract

\endtopmatter

\document
\heading {\S 1. Introduction and results} \endheading

For more than hundred years Schubert decomposition of complete and partial 
flag manifolds and Schubert calculus has been an intertwining 
point of algebra, geometry, representation theory and combinatorics.
In this paper we are concerned with the further refinement of Schubert
decomposition theory.

 Our far and ambitious goal is to describe the topology of intersection 
of two arbitrary Schubert cells. This problem seems to be very important
in connection with representation theory, in particular, with calculation 

of Kazhdan--Lusztig  polynomials, which is known to be a very hard problem.

As a minor step in this direction we calculate in this paper 
the number of connected components in the intersection of two open opposite 
Schubert cells in the space of real $n$-dimensional flags. 
This question was raised in \cite {A},
and later, in connection with criteria of total positivity, in \cite {SS}.
Other topological characteristics of pairwise intersections of Schubert cells
(not neccessary opposite) are considered in our previous papers
\cite{SV, SSV1, SSV2}. 

Despite the fact that in this paper we consider only real algebraic varieties,
we hope that our approach can give some information about complex geometry
as well. The reasons for this hope are as follows.

a) The calculation method  is based on the beautiful and very
deep (``very algebraic'') results by Berenstein, Fomin and Zelevinsky 
\cite{BFZ}
on the Lusztig parametrization of the unipotent lower triangular
matrices, which form an open cell in the flag manifold. It looks like
the same type of parametrization and cell decompositions 
might be used also for calculation
of certain topological characteristics (e.g., the fundamental group)
in the complex case.

b) Any real algebraic manifold satisfies the so-called Smith inequality:
the sum of its Betti numbers for homology with $\bF_2$-coefficients
is less or equal than the sum of Betti numbers of its complexification
with $\Bbb Z$-coefficients. By the results of one of the authors (see
\cite{S}), the intersection of opposite open cells in the space of complete
flags enjoys the so-called  M-property, i.e., the Smith
inequality becomes an equality. It is another indirect hint that 
real and complex geometries of this intersection are related very closely.  

Recall that a real flag $f$ is a sequence of $n$ linear 
subspaces $\{f^1\subset f^2\subset \dots\subset f^{n-1}\subset f^n=\bR^n\}$,
where $\dim f^k=k$. We denote the set of all real flags by $F_n$. 
In an obvious way one can identify $F_n$ with the set $SL_n(\bR)/B$, 
where $B$ is the Borel subgroup of all upper triangular matrices.
Indeed, fix coordinates in $\bR^n$. 
Then any nondegenerate matrix produces a flag whose $i$-dimensional
subspace is generated by the first $i$ columns of the matrix. 
 
Two flags $f$ and $g$  are called transversal if
all linear subspaces of $f$ are transversal to all subspaces
of $g$, i.e.,  for any pair $1\ls i,j\ls n$ 
$$
\dim (f^i\cap g^j)=\left\{
\aligned 
i+j-n, &\quad\text{ if } i+j\gs n,\\
  0, &\quad\text{ otherwise. } 
\endaligned
\right.
$$

Let us choose a pair of transversal flags $f$ and $g$.
The main object of our interest is the set $U_{f,g}^n$ 
of all flags transversal to both $f$ and $g$.
Evidently, $U_{f,g}^n$ is the intersection of open maximal
Schubert cells in Schubert decompositions of the
variety of complete real flags $SL_n(\bR)/B$ 
with respect to the flags $f$ and $g$. Since the group $SL_n(\bR)$ acts 
transitively on pairs of transversal flags, all $U_{f,g}^n$ are diffeomorphic.
We thus omit the indices $f$ and $g$ and write just $U^n$.

Let us now use the isomorphism $F_n\cong SL_n(\bR)/B$ to 
introduce coordinates in $F_n$ (and, in particular, in $U^n$).
Denote by $\{e_1,\dots, e_n\}$ the standard basis in $\bR^n$ and take for
$f$ the flag spanned by  $\{e_1\}, \{e_1,
e_2\},\dots$, and for $g$ the flag spanned by $\{e_n\},\{e_n,e_{n-1}\},\dots$.
Then the open Schubert cell of all flags transversal to $g$ is 
identified naturally with the set of all lower triangular unipotent 
matrices; namely, the $i$th subspace of a variable flag is spanned by the first
$i$ columns of the matrix.

In order to make our notation consistent with that of \cite{BFZ} we
transpose lower triangular matrices and parametrize
flags transversal to $g$ with the set $N^n$ of all upper 
triangular unipotent matrices. 

Flags transversal to both $f$ and $g$ are described in the following way.
Let  $\l$ be an arbitrary subset of $\{1,2,\dots,n\}$,
and $D_\l$ be the determinant of the submatrix formed by the
first $|\l|$ rows and the columns numbered by the elements of $\l$;
we write $D_k$ instead of $D_{\{n-k+1,n-k+2,\dots,n\}}$. 
The set of all flags transversal
to both $f$ and $g$ corresponds to the set of all upper triangular 
matrices satisfying the additional restriction $D_i\neq 0$ for all 
$i=1,\dots,n-1$ (see \cite{SV} for details).
Let us denote by $\Delta_i$ the divisor $\{D_i=0\}\subset N^n$, and
by $\Delta^n$ the union $\cup_{i=1}^{n-1} \Delta_i$.
We thus see that $U^n$ is diffeomorphic to the complement
$N^n\setminus\Delta^n$.

The main result of the present paper is a reduction of the topological problem
of finding the number of connected components in the intersection of two
open opposite Schubert cells to a purely combinatorial problem of enumerating
the orbits of a certain linear group action in the vector space of 
upper triangular matrices with $\bF_2$-entries. 

This combinatorial reduction consists of two steps.
The first one seems very natural and was obtained independently by K.Rietsch 
in a more general situation (for the case of A- and D-series).
It leads to a combinatorial problem 
of enumerating the connected components of a certain graph,
see \cite{R} and \S 2. The second step allows to 
decrease significantly the cardinality of the graph considered  and, finally,
leads to a $\bF_2$-linear group action.  

The resulting reduction of the initial problem to the group action can be
generalized in a natural way to the case when the initial flags $f$ and $g$ 
are not in general position. Presumably, the number of orbits of this 
generalized group action coincides with the number of connected components 
in $U^n_{f,g}$ for any
relative position of $f$ and $g$. Strangely enough, the first step of 
the reduction does not generalize directly to the case of an arbitrary 
pair $(f,g)$. Apparently, 
there exists a deeper connection between the initial problem and the group
action than the one known to the authors for the present moment, see
Final Remarks. 

 Let us formulate the main result of our paper.
Consider the vector space $N^n(\bZ_2)$ of upper triangular matrices with
$\bF_2$-valued entries. We define the group $\G_n$ as the subgroup of
$GL(N^n(\bF_2))$ generated by $\bZ_2$-linear transformations $g_{ij}$, 
$1\ls i\ls j\ls n-1$. The generator $g_{ij}$ acts on a matrix 
$M\in N^n(\bF_2)$ as follows. Let  $M^{ij}$ denote the $2\times 2$ minor
of $M$ formed by  rows $i$ and $i+1$ and columns $j$ and $j+1$
(or its upper triangle in case $i=j$). Then $g_{ij}$
applied to $M$  changes  $M^{ij}$ by adding to each entry of $M^{ij}$ the 
$\bF_2$-valued trace of $M^{ij}$, and does not change all the other entries
of $M$. For example, if $i<j$ and
$
M^{ij}=\pmatrix a & b \cr
        c & d
\endpmatrix,
$ 
then $g_{ij}$ changes $M^{ij}$ as follows:
$$\pmatrix a & b \cr
        c & d 
\endpmatrix 
\mapsto 
\pmatrix d & b+a+d \cr
        c+a+d & a 
\endpmatrix.
$$

Observe that each $g_{ij}$ is an involution on $N^n(\bF_2)$.
Some other properties of $\G_n$ are given in Final Remarks.
Apparently, an explicit description of $\G_n$ can be derived
from the above mentioned properties and the classification contained
in \cite{J1}. The most essential result relating
the properties of $G_n$ to our initial question is as follows. 

\proclaim{Main Theorem}
The number $\s_n$ of connected components in $U^n$ 
coincides with the number of $\G_{n-1}$-orbits in $N^{n-1}(\bZ_2)$.
\endproclaim

The main conjecture, which is based on our computer
experiments and a detailed study of the $\G_n$-action, is as follows.

\proclaim{Main Conjecture} 
The number $\sharp_n$ of connected components in $U^n$ equals
$3\times 2^{n-1}$ for all $n> 5$. 
\endproclaim

Cases $n=3,4 \text{ or }5$ are exceptional, with $\sharp_3=6$, $\sharp_4=20$, 
$\sharp_5=52$.

The structure of the paper is as follows. In \S 2 we recall some basic
constructions
of \cite {L} and \cite{BFZ} and get the first combinatorial reduction of the
problem. Section \S 3 is central in the paper. Here we describe the second
combinatorial
reduction and prove that it is consistent.

We are very grateful to A.~Zelevinsky for a number of stimulating
discussions and encouragement, and to K.~Rietsch who (after the talk 
on the topic given by B.~Shapiro at the Northeastern University)
provided us with a copy of her unpublished paper \cite{R}.
B.~Shapiro and A.~Vainshtein express their gratitude to Volkswagen-Stiftung 
for the financial support of their stay at MFO in Oberwolfach 
(under the program ``Research in Pairs''),
where an essential progress in the proof of the above main conjecture was
achieved.

\heading {\S 2. Chamber ansatz and the first combinatorial reduction} 
\endheading

\subheading{2.1}
Significant recent papers \cite {L} and \cite {BFZ} have shown the 
importance of the following parametrization of the set $N_{>0}^n$ of 
totally positive real upper triangular $n\times n$
matrices. In particular, it implies a number of 
interesting results on canonical bases of quantum groups, as well as new 
criteria of total positivity. Namely, denote
by $w_0^n$ the permutation of the maximal length in the symmetric group $S_n$ 
and fix some reduced decomposition $w_0^n=s_{i_1}\ldots s_{i_m}$,
where $s_{i_j}$ stands for the transposition $(i_j,i_j+1)$, $1\ls i_j\ls n-1$,
$m=n(n-1)/2$. Following  G.~Lusztig, one can factorize any matrix 
$\Cal M\in N_{>0}^n$ into the product
$$
\Cal M=(\one + t_1\cdot E_{s_{i_1}})(\one + t_2\cdot E_{s_{i_2}})\times
\cdots\times(\one + t_m\cdot E_{s_{i_m}}), \tag1
$$ 
where $E_{s_{i_l}}$ is the matrix with the only nonvanishing entry 
$e_{s_{i_l},s_{i_l}+1}=1$, and all $t_l$, $l=1,\dots,m$, are positive.

The right hand side of the above formula yields a well-defined mapping
$L_\si\:\BbbR^m\to N^n$, where
$\si$ stands for the reduced decomposition of $w_0^n$ chosen above.

A.~Berenstein, S.~Fomin and A.~Zelevinsky have derived the formulas for the
parameters $\{t_l\}$ in terms of the matrix entries of $\Cal M$.
Using their formulas, or directly, one can get the following expressions for
$D_i(\Cal M)$, which hold true for any reduced decomposition $\si$:
$$
\align
&D_1(\Cal M)=t_1\cdot\ldots\cdot t_n,\\
&D_2(\Cal M)=t_2\cdot\ldots\cdot t_n\cdot t_{n+1}\cdot\ldots\cdot
t_{2n-1},\\ 
&\qquad\dots\\
&D_{n-1}(\Cal M)=\prod_{i=1}^{n-1} t_{kn-\frac{k(k+1)}{2}}.\endalign 
$$

Observe that each $D_i$ is a monomial in $\{t_i\}$. 
Therefore, if $t_i\ne 0$ for all $i=1,\dots,m$, then all 
the $D_i$ do not vanish. Thus, the image $\Cal U_\si$ of
the map $L_\sigma\: (\BbbR\setminus 0)^m\to N^n$ for any reduced 
decomposition $\sigma$ of $w_0^n$ lies in $U^n$. In fact, the explicit
expressions for $\{t_l\}$ imply the following proposition.

\proclaim{Proposition} $L_\si$ is a 
diffeomorphism of $(\BbbR\setminus 0)^m$ onto its image $\Cal U_\si
\subset U^n$.
\endproclaim

Denote by $C_\sigma$ the complement to $\Cal U_\si$ in $U^n$, and
let $\Sigma^n$ denote the set of all reduced decompositions of 
$w_0^n$.

\proclaim{2.2. Lemma} The codimension of  
$\bigcap_{\sigma\in\Sigma^n}C_\sigma$ in $U^n$ is at least $2$. 
\endproclaim

\demo{Proof}
To prove the claim it suffices to find two reduced decompositions $\sigma$ 
and $\tau$ such that the codimension of $C_\sigma\cap C_\tau$ is already
at least 2. It follows from \cite {BFZ} that $C_\sigma$ for any $\sigma\in
\Sigma^n$ is represented as the union of irreducible divisors 
given by equations $M_\l=0$, where the index $\l$ runs over a certain subset
$\L_\sigma$ of increasing subsequences of $\{1,2,\dots,n\}$. The explicit 
expression for $M_\l$ is given by the so-called Chamber Ansatz of \cite {BFZ}.
Namely, let $\Cal M$ be an $n\times n$ matrix and $[\Cal M]_+$ denote the
last factor $\Cal M_2$ in the Gaussian LDU-decomposition $\Cal M=
\Cal M_1^T\Cal D\Cal M_2$, where $\Cal M_1,\Cal M_2\in N^n$ and $\Cal D$ is
diagonal. Then $M_\l(\Cal M)=D_\l(\Cal N)$, where $\Cal M=[w_0^n\Cal N^T]_+$.

Let us fix 
$\sigma=s_1s_2\dots s_{n-1}s_1s_2\dots s_{n-2}\dots s_1s_2s_1$ and
$\tau=s_{n-1}s_{n-2}\dots s_1s_{n-1}s_{n-2}\dots 
\allowmathbreak s_2\dots s_{n-1}s_{n-2}s_{n-1}$. 
It follows from the Chamber Ansatz that
$$\align
\L_\sigma&=\{\{r,r+1,\dots,q\}\: 1<r<q<n\},\\
\L_\tau&=\{\{1,2,\dots,s\}\cup\{p,p+1,\dots,n-1\}\: 0<s<p<n\},
\endalign
$$
and thus $\L_\sigma\cap\L_\tau=\varnothing$.

Consider a matrix $\Cal M$ of the form 
$$
\Cal M=\pmatrix
1 & * & * & * & \dots & * & 1 \cr
0 & 1 & * & * & \dots & * & 0 \cr
0 & 0 & 1 & * & \dots & * & 0 \cr
\vdots & & &  & \ddots& &\vdots \cr
0 & 0 & 0 & 0 & \dots & 1 & 0 \cr
0 & 0 & 0 & 0 & \dots & 0 & 1 
\endpmatrix.
$$

It is easy to check that for a generic $\Cal M$ one has
$M_\l(\Cal M)=0$ for any $\l\in\L_\sigma$, but
$M_\l(\Cal M)\ne 0$ for $\l\in\L_\tau$. 
Since the divisors $M_\l=0$ are irreducible for all 
$\l\in\L_\sigma\cup\L_\tau$, we thus get
$\codim C_\sigma\cap C_\tau\gs 2$.
\qed
\enddemo

\proclaim{2.3. Corollary} The number $\s_n$ is equal to the number of 
connected components of $\cup_{\sigma\in\Sigma^n}\Cal U_\si$.
\endproclaim

\demo{Proof} The number of connected components does not change if we delete a
subset of codimension $\gs2$.
\qed
\enddemo

\subheading {2.4. Transition rules}
The coefficients $\{t_i\}$ of factorization (1) depend on the choice of 
the reduced decomposition $\si$. The corresponding transition formulas,
borrowed from \cite{BFZ}, are given below. All reduced
decompositions can be obtained from each other by
a sequence of the so-called 2- and 3-moves. A 2-{\it move\/} is the
interchange of neighboring $s_{i_k}$ and $s_{i_{k+1}}$ in the decomposition 
under the assumption that they commute,
and a 3-{\it move\/} is the substitution of $s_{j+1}s_js_{j+1}$ instead of
$s_js_{j+1}s_{j}$.

The transition formulas for the 2-move at the positions $(i,i+1)$ are as 
follows:
$$
\aligned &t'_i=t_{i+1},\\
&t'_{i+1}=t_i.\endaligned\tag2
$$

The transition formulas for the 3-move at the positions $(i,i+1,i+2)$ are
more complicated:
$$\aligned &t'_i=\frac{t_{i+1}t_{i+2}}{t_i+t_{i+2}},\\
&t'_{i+1}=t_i+t_{i+2},\\
&t'_{i+2}=\frac{t_{i+1}t_i}{t_i+t_{i+2}}.\endaligned \tag3
$$

\subheading{2.5. Sign transition rules}
Let $\kap_i$ stand for the sign of $t_i$, that is, $\kap_i\in\{+,-\}$ (we
assume that $t_i$ does not vanish). Below we formulate transition rules 
for $\{\kap_i\}$ under any 2- or 3-move.

For any 2-move, the new signs are defined uniquely from (2):
$$
\aligned &\kap'_i=\kap_{i+1},\\
&\kap'_{i+1}=\kap_i.\endaligned \tag4
$$

For 3-moves, there are two different possibilities:

1) in the following cases the new signs are defined uniquely from (3):
$$
\matrix
(\kap_i' &\kap_{i+1}'&\kap_{i+2}')&\leftrightarrow  &
(\kap_i  &\kap_{i+1} &\kap_{i+2})\cr
(+       &+          &\quad+)     &\leftrightarrow  &
(+       &+          &\quad+)\cr
(-       &-          &\quad-)     &\leftrightarrow  &
(-       &-          &\quad-)\cr
(+       &-          &\quad+)     &\leftrightarrow  &
(-       &+          &\quad-)\cr
(-       &+          &\quad-)     &\leftrightarrow  &
(+       &-          &\quad+)
\endmatrix;\tag5
$$

2) for all the other cases, the changes of signs are given in the following
table:
$$
\matrix
(\kap_i' &\kap_{i+1}'&\kap_{i+2}')&\leftrightarrow  &
(\kap_i  &\kap_{i+1} &\kap_{i+2})\cr
         &           &            &                 &
(+       &-          &\quad-)\cr
(+       &+          &\quad-)     &\leftrightarrow  &
         &\text{or}  &\cr
         &           &            &                 &
(-       &+          &\quad+);\cr
         &           &            &                 &
(+       &-          &\quad-)\cr
(-       &-          &\quad+)     &\leftrightarrow  &
         &\text{or}&\cr
         &           &            &                 &
(-       &+          &\quad+);\cr
         &           &            &                 &
(+       &+          &\quad-)\cr
(+       &-          &\quad-)     &\leftrightarrow  &
         &\text{or}  &\cr
         &           &            &                 &
(-       &-          &\quad+);\cr
         &           &            &                 &
(+       &+          &\quad-)\cr
(-       &+          &\quad+)     &\leftrightarrow  &
         &\text{or}&\cr
         &           &            &                 &
(-       &-          &\quad+).\cr
\endmatrix\tag6
$$

\proclaim{2.6. Lemma} All the changes of signs listed in {\rm (6)}
can be realized. \endproclaim

\demo{Proof}
It suffices to check only the case
$$
\matrix
         &           &            &                 &
(+       &-          &\quad-)\cr
(+       &+          &\quad-)     &\leftrightarrow  &
         &\text{or}  &\cr
         &           &            &                 &
(-       &+          &\quad+).\cr
\endmatrix
$$
Let us prove, for example, that the sign sequence $(++-)$ can be
transformed into both $(+--)$ and $(-++)$.
Suppose that $t_i+t_{i+2} > 0$. Then one has $t_i>0$, $t_{i+1}>0$, 
$t_{i+2}<0$, and $t_i+t_{i+2} > 0$. Transition formulas (3) imply 
$t'_i=\frac{t_{i+1}t_{i+2}}{t_i+t_{i+2}}<0$. In the same way one gets from
(3) $t'_{i+1}>0$ and  $t'_{i+2}>0$, and thus we get the transition
$(++-)\leftrightarrow (-++)$. Since $t_i$ and $t_{i+2}$ have opposite 
signs, we can deform the initial parameters $t_i$, $t_{i+1}$, $t_{i+2}$ 
(in a nonvanishing way) into $\widetilde t_i$, $\widetilde t_{i+1}$,
$\widetilde t_{i+2}$ such that $\widetilde t_i+\widetilde t_{i+2} < 0$. 
This gives us the second transition $(++-)\leftrightarrow (+--)$.

In the same way one can show that all the other transitions listed 
in (6) can be realized as well.
\qed \enddemo

\subheading{2.7. The graph $G^n$ of reduced decompositions modulo 2-moves}
Consider the set $G^n$ of all reduced decompositions of $w_0^n$ modulo 
2-moves. Let us present $G^n$ as a graph (which we also denote by $G^n$). 
The vertices of $G^n$ are the equivalence  classes of reduced decompositions 
modulo 2-moves. Two vertices are connected by an edge if there exists a 
pair of representing reduced decompositions such that some 3-move sends 
one of them to the other one.

Another well-known interpretation of $G^n$ is the set of topologically
different arrangements of $n+1$ pseudolines in $\Bbb R^2$. The set $G^n$ 
was studied, e.g., in \cite {OM, pp.247--280} and \cite {Kn, pp.29--40}, 
but to the  best of the authors knowledge, even
the cardinality of $G^n$ is still unknown for large $n$.

Observe that if two reduced decompositions $\si_1$ and $\si_2$ are 
equivalent modulo 2-moves, then the corresponding sets $\Cal U_{\si_1}$
and $\Cal U_{\si_2}$ coincide, by (2). Let us fix an arbitrary
representative $\sigma(v)$ in each equivalence class $v$. Then it
follows from Corollary 2.3 that $\s_n$ is equal to the number of connected 
components  of $\cup_{v\in G^n}\Cal U_{\si(v)}$.

\subheading{2.8. The ``large'' graph $\widetilde G^n$ as a covering of $G^n$}
Let us now construct a certain covering $\widetilde G^n$ of $G^n$. Namely, a
vertex of $\widetilde G^n$ is a pair consisting of a vertex $v$ of $G^n$ and 
a set of $m=n(n-1)/2$ signs (that is, $+$'s or $-$'s) interpreted as values
of the variables $\{\kap_i\}$ for the decomposition $\si(v)$. Thus, the
fiber over any vertex of $G^n$ contains exactly $2^m$ vertices of
$\widetilde G^n$. Two vertices of $\widetilde G^n$
are adjacent if

1) their projections in $G^n$ are adjacent;

2) the corresponding variables $\kap_i,\kap_{i+1},\kap_{i+2}$ and
$\kap'_i,\kap'_{i+1},\kap'_{i+2}$ satisfy  relations (5) or (6).

Let us denote by $\pi\:\widetilde G^n\to G^n$ the natural projection.

\proclaim {2.9. Theorem (first combinatorial reduction)} The number $\s_n$
of connected components of $U^n$ is equal to that of $\widetilde G^n$. 
\endproclaim

\demo{Proof} A vertex $(v,\kap)$, $\kap=\{\kap_1,\dots,\kap_m\}$, of
$\widetilde G^n$ can be identified with the image $\Cal U_{\si(v)}^\kap$ of 
the set $\bR_{\kap_1}\times\cdots\times\bR_{\kap_m}$ under the mapping
$L_{\si(v)}$. By 2.7, $\s_n$ is just the number of connected components
in the union of the above images.

Suppose that $(v,\kap)$ and $(v',\kap')$ are adjacent in $\widetilde G^n$.
Then there exist $\hat\si\in v$ and $\hat\si'\in v'$ such that $\hat\si'$
is obtained from $\hat\si$ by a 3-move and the corresponding variables
$\hat\kap$ and $\hat\kap'$ (obtained from $\kap$ and $\kap'$, respectively,
via (4)) satisfy (5) or (6). Similarly to the proof of Lemma 2.6, we see
that the intersection of $\Cal U_{\hat\si}^{\hat\kap}$ and
$\Cal U_{\hat\si'}^{\hat\kap'}$ is nonvoid. However, 
$\Cal U_{\hat\si}^{\hat\kap}=\Cal U_{\si(v)}^\kap$ and
$\Cal U_{\hat\si'}^{\hat\kap'}=\Cal U_{\si(v')}^{\kap'}$. Thus,
$\Cal U_{\si(v)}^\kap\cap\Cal U_{\si(v')}^{\kap'}\ne\varnothing$,
and hence $\Cal U_{\si(v)}^\kap$ and $\Cal U_{\si(v')}^{\kap'}$
belong to the same connected component of $U^n$.

Suppose now that 
$\Cal U_{\si(v)}^\kap\cap\Cal U_{\si(v')}^{\kap'}\ne\varnothing$.
Then, by 2.1, there exists a sequence $\{\si_1,\si^1,\si_2,\si^2,\dots,
\si_k,\si^k\}$ such that $\si_1=\si(v)$, $\si^k=\si(v')$, $\si_j$ and $\si^j$
are equivalent modulo 2-moves, and $\si_j$ is obtained from $\si^{j-1}$
by a 3-move. Denote by $v''$ the equivalence class of $\si_2$ modulo 2-moves
and fix a generic matrix 
$\Cal M\in \Cal U_{\si(v)}^\kap\cap\Cal U_{\si(v')}^{\kap'}$. 
The matrix $\Cal M$
and relations (4)--(6)  define uniquely $\kap''$ such that the vertices
$(v,\kap)$ and $(v'',\kap'')$ are adjacent in $\widetilde G^n$. Proceeding
in the same way, we get a path from $(v,\kap)$ to some $(v',\hat\kap')$ in
$\widetilde G^n$ such that
$\Cal M\in \Cal U_{\si(v)}^\kap\cap\Cal U_{\si(v')}^{\hat\kap'}$. However, the
sets $\Cal U_{\si(v')}^{\kap'}$ and $\Cal U_{\si(v')}^{\hat\kap'}$ for 
$\kap'\ne\hat\kap'$ are disjoint, and thus $\kap'=\hat\kap'$. Therefore, 
$(v,\kap)$ and $(v',\kap')$ lie in the same connected component of 
$\widetilde G^n$. 
\qed
\enddemo

The above theorem reduces our initial problem to a purely combinatorial
setup. However, since
the number of vertices in $\widetilde G^n$ equals $2^{n(n-1)/2}\cdot|G^n|$, 
a direct solution of this combinatorial problem is hardly possible.

\heading \S 3. Reduction to the fiber of $\pi\:\widetilde G^n\to G^n$. 
\endheading

Theorem 3.2 below reduces our initial problem to the 
calculation of the number of connected components in a smaller graph 
$\Gamma^n$ with ``only'' $2^{n(n-1)/2}$ vertices.

\subheading{ 3.1. Construction of the graph $\Gamma^n$}
Let us fix the vertex $v_0^n\in G^n$ that corresponds to the equivalence class 
of the reduced decomposition $s_1s_2\dots s_{n-1}s_2s_2\dots s_{n-2}\dots
s_1s_2s_1$. 

\vskip 10pt
\centerline{\hbox{\hskip 0.2cm \epsfysize=3cm\epsfbox{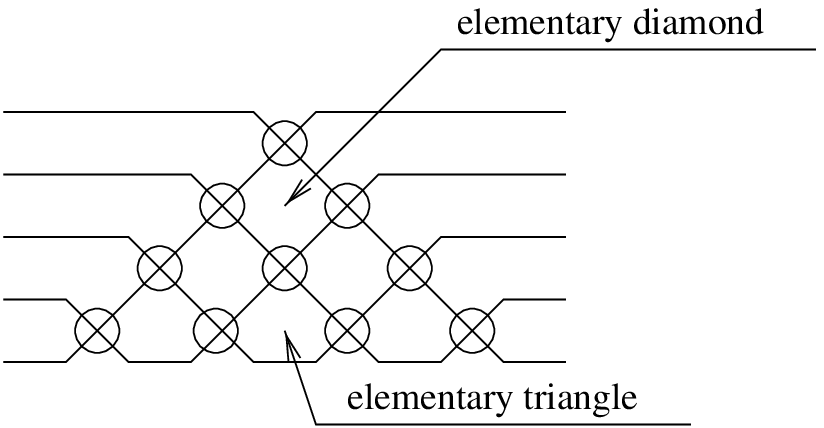}}}
\midspace{1mm} \caption{Fig.1. The special vertex $v_0^n$} \vskip 5pt

We consider $v_0^n$ as an arrangement of $n$ pseudolines (see Fig.~1). 
The intersection point of two pseudolines is called a {\it node\/};
thus, $v_0^n$ has $m=n(n-1)/2$ nodes that lie on $n-1$ {\it horizontal
levels\/} (the $i$th level contains $n-i$ nodes).  
Define the {\it vertex set of\/} $\Gamma^n$ as the set $\{+,-\}^m$ of all 
the choices of signs at all the nodes of $v_0^n$.  Obviously, this set is 
one of the fibers of the projection
$\pi\: \widetilde G^n\to G^n$, namely, the inverse image of $v_0^n$.

Let us now define the {\it edge set of\/} $\Gamma^n$. The faces of the 
arrangement $v_0^n$ are of two types: triangles and diamonds (see Fig.~1). 
Let us call them {\it elementary regions}. Observe that each elementary
region has exactly two {\it horizontal\/} nodes (those lying on the same
horizontal level). The number of this level is called the {\it height\/}
of the elementary region.

For any elementary region $A$ we define the 
involution $I_A$ on the vertex set of $\Gamma^n$: $I_A$ reverses 
the signs at all the nodes of $A$. An involution $I_A$ is called {\it 
admissible\/} for a vertex $\kap\in\Gamma^n$ if the horizontal nodes of
$A$ have the opposite signs in $\kap$. Two vertices $\sg_1$ and $\sg_2$
of $\Gamma^n$ are connected by an edge if  and only if there exists an
elementary region $A$ such that $\sg_1=I_A(\sg_2)$ and $I_A$
is admissible for $\kap_1$ (and thus, for $\kap_2$). 

\proclaim{3.2. Theorem (second combinatorial reduction)} 
The number of connected components of $\widetilde G^n$ is equal to that 
of $\Gamma^n$. \endproclaim

The proof of this main result of the paper splits up into two major parts 
and is postponed till the end of this section.

\proclaim{3.3. Proposition} Let $\kap$ and $\kap^*$ belong 
to the same connected component of $\Gamma^n${\rm,} then $(v_0^n,\kap)$
and $(v_0^n,\kap^*)$ can be connected by a path in $\widetilde G^n$. 
\endproclaim

\demo{Proof} Without loss of generality one can assume that $\sg$ and 
$\sg^*$ are the endpoints of an edge in $\Gamma^n$. Below we construct certain
paths in $\widetilde G^n$ connecting $(v_0^n,\kap)$ and $(v_0^n,\kap^*)$.

Let $A$ be the elementary region that defines the edge $[\sg,\sg^*]$. 
We proceed by induction on its height $h(A)$.

Let $h(A)=1$, then $A$ is an elementary triangle. Denote
the signs of its nodes in $\sg$ by $\kap_i$, $\kap_{i+1}$, $\kap_{i+2}$.
Since $[\sg,\sg^*]\in\Gamma^n$, we have $\kap_i=\overline\kap_{i+2}$ 
(hereinafter
$\overline\kap$ stands for the sign opposite to $\kap$). We thus can use 
the rules (6) twice to get first $\kap'_i=\kap_i$, 
$\kap'_{i+1}=\overline\kap_{i+1}$, $\kap'_{i+2}=\kap_{i+2}$, and then
$\kap''_i=\overline\kap'_i=\overline\kap_i$, $\kap''_{i+1}=\kap''_{i+1}=
\overline\kap_{i+1}$, $\kap''_{i+1}=\overline\kap'_{i+2}=
\overline\kap_{i+2}$. Hence, the path between $(v_0^n,\kap)$ and 
$(v_0^n,\kap^*)$ $\widetilde G^n$
in this case consists of two edges (see Fig.~2 for an illustration).
Observe that all the nodes involved in the above transitions lie on the
two pseudolines containing the left side and the base of $A$.

\vskip 10pt
\centerline{\hbox{\hskip 0.2cm \epsfxsize=10truecm\epsfbox{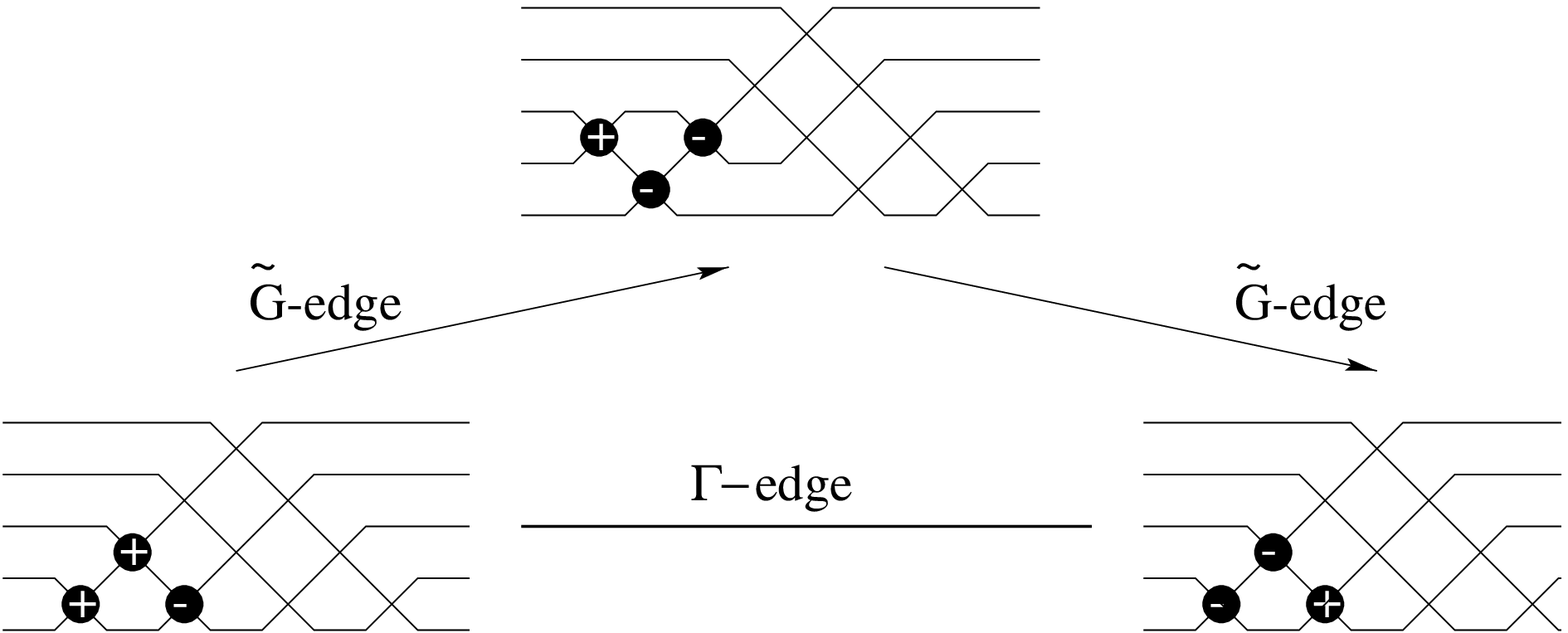}}}
\midspace{1mm}
\caption{Fig.2. The sign flip for the nodes of an elementary triangle} 
\vskip 5pt

Assume now that the assertion holds whenever the height of the elementary
region $A$ defining the edge $[\sg,\sg^*]$ is at most $k-1$; moreover, 
assume
that all the nodes involved in the neccessary transitions lie on the
two pseudolines containing the upper left and the lower right sides of $A$.
 
Let now $h(A)=k>1$, and hence $A$ is an elementary diamond.
Consider the sequence $\{A_1,\dots,A_k\}$ of elementary regions
defined as follows:

$A_1=A$;

the upper node of $A_i$ is the leftmost node of $A_{i-1}$, $2\ls i
\ls k$;

the rightmost node of $A_i$ is the lower node of $A_{i-1}$, $2\ls i
\ls k$.

Evidently, $A_k$ is an elementary triangle. Applying (5) or (6) to 
the signs of its nodes, we can pass
to a neighboring vertex $(v_1^n,\kap')$ of $\widetilde G^n$ (see Fig.~3).

\vskip 10pt
\centerline{\hbox{\hskip 0.2cm \epsfxsize=13truecm\epsfbox{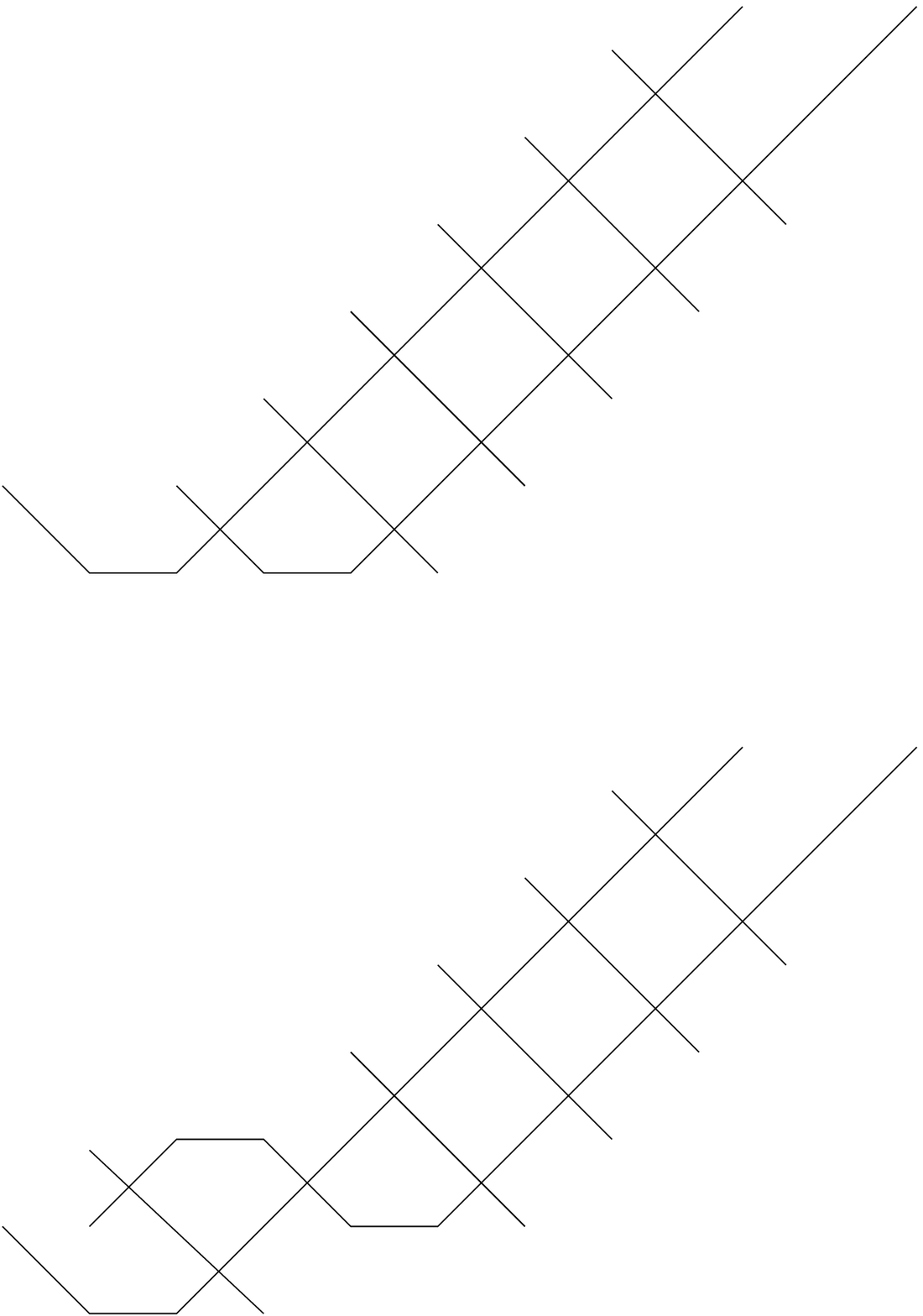}}}
\midspace{1mm}
\caption{Fig.3. To the proof of Proposition 3.3}
\vskip 5pt

Let us now delete from the arrangement $v_1^n$ the pseudoline that defines
the right side of the triangle $A_k$ in $v_0^n$. We thus get an
arrangement of $n-1$ pseudolines, which is evidently isomorphic to 
$v_0^{n-1}$,  and a sign distribution $\tilde\kap'$ on it, which corresponds
naturally to $\kap'$. However, the height of $A$ in this arrangement
equals $k-1$, and thus, there exists a path from $(v_0^{n-1},\tilde\kap')$ to 
$(v_0^{n-1},\tilde\kap'')$ in the graph $\widetilde G^{n-1}$, where
$\tilde\kap''=I_A(\tilde\kap')$ in $v_0^{n-1}$. Since
this path involves only nodes that lie on the two pseudolines defining the
upper left and the lower right sides of $A$, this path is lifted 
naturally to a path in $\widetilde G^n$ that connects $(v_1^n,\kap')$ with
$(v_1^n,\kap'')$ such that $\kap''=\overline\kap'$ at the nodes of $A$
and $\kap''=\kap'$ otherwise. To complete the induction step it suffices
to apply in the opposite direction the transition that leads from 
$(v_0^n,\kap)$ to $(v_1^n,\kap')$ and to note that the nodes involved in this
transition lie on the above described pseudolines (see Fig.~4 for an
illustration).
\qed
\enddemo

\vskip 10pt
\centerline{\hbox{\hskip 0.2cm \epsfxsize=10truecm\epsfbox{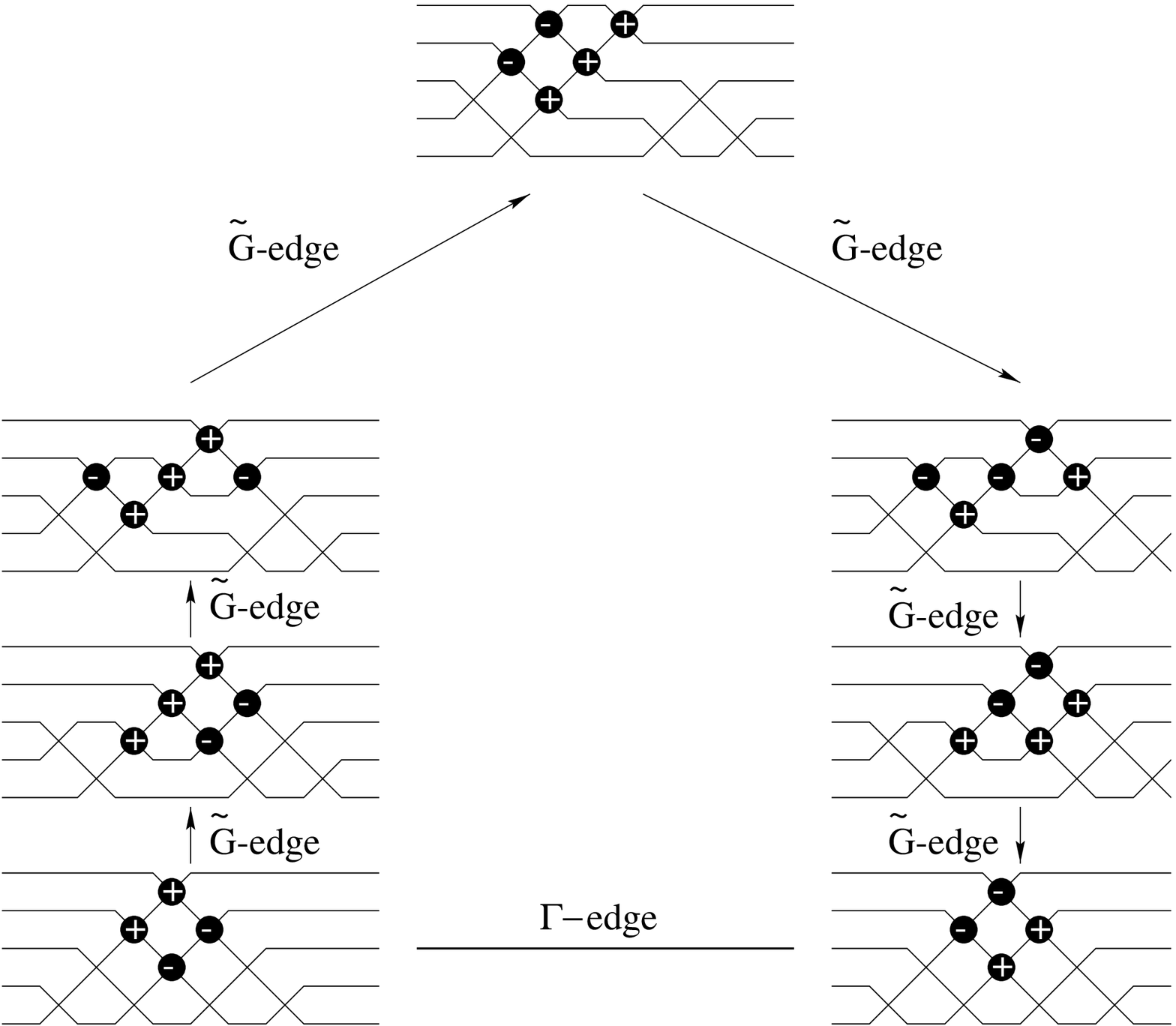}}}
\midspace{1mm}
\caption{Fig.4. The sign flip for the nodes of an elementary diamond}
\vskip 5pt

\subheading{3.4. Outline of the further strategy} Our goal is to prove the 
inverse of Proposition~3.3, namely, that if two vertices of $\Gamma^n$ are
the endpoints of a path $p$ in $\widetilde G^n$, then they lie in the same 
connected component of $\Gamma^n$ (see Proposition 3.16). Consider the
projection $\pi(p)$ of such a path $p$; evidently, it is a closed path in
$G^n$. We next prove the following property of the paths $p$ in
$\widetilde G^n$ whose projection $\pi(p)$ is a closed path in $G^n$: for any 
$p$ satisfying the above condition there exists a sequence of paths in
$\widetilde G^n$ such that each path in the sequence has
a simple structure (canonical lifts and special lifts, see 3.5 and 3.6), and
the concatenation of all the paths in the sequence is a path between 
the endpoints of $p$ (see Lemma 3.7). Therefore, it remains to prove that 
the inverse of Proposition~3.3 is true for special lifts (Corollary 3.10) 
and for canonical
lifts of closed paths (Corollary 3.15). To obtain the first of these results
we establish certain relations between two fibers $\Gamma^n(u)$ and
$\Gamma^n(v)$ over adjacent vertices $u$ and $v$ of $G^n$ (see 3.8 and 
Lemma~3.9). To get the second one, we find a basis of the space of closed
paths in $G^n$ (see 3.11 and Lemma~3.14) and prove the corresponding statement
for all the elements of the basis (see Lemmas~3.12 and~3.13).

\subheading{3.5. Canonical lifts} Consider an arbitrary edge $e=[u^1,u^2]$ 
in $G^n$ and an arbitrary vertex $\tilde u^1=(u^1,\kap^1)\in\pi^{-1}(u^1)$ 
of $\widetilde G^n$. The sign transition defined by the 3-move
$e$ is governed either by (5), or by (6). In the first case, there exists 
a unique vertex $\tilde u^2=(u^2,\kap^2) \in\pi^{-1}(u^2)$
such that $[\tilde u^1,\tilde u^2]$ is an edge of $\widetilde G^n$. We then say
that  $[\tilde u^1,\tilde u^2]$ is the {\it canonical
$\tilde u^1$-lift\/} of $e$ (see Fig.~5). 

In the second case, there exist two
vertices $\tilde u^2,\tilde v^2\in\pi^{-1}(u^2)$ such that both 
$[\tilde u^1,\tilde u^2]$ and $[\tilde u^1,\tilde v^2]$ are edges in 
$\widetilde G^n$; besides,  there 
exists also $\tilde v^1\in\pi^{-1}(u^1)$ such that both 
$[\tilde v^1,\tilde u^2]$ and $[\tilde v^1,\tilde v^2]$ are edges in 
$\widetilde G^n$. In this case, the  canonical
$\tilde u^1$-lift of $e$ is the edge that preserves the middle element 
$\kap_{i+1}$ of the corresponding triple at $\tilde u^1=(u^1,\kap^1)$ 
(see (6)). Observe that in this case the canonical $\tilde v^1$-lift of $e$ 
is the only edge among  $[\tilde v^1,\tilde u^2]$ and $[\tilde v^1,\tilde v^2]$
that does not touch the  canonical $\tilde u^1$-lift of $e$. The canonical 
$\tilde u^1$-  and $\tilde v^1$-lifts of $e$ are called {\it twins} (see
Fig.~5).

\vskip 20pt
\centerline{\hbox{\hskip 0.2cm \epsfxsize=15truecm\epsfbox{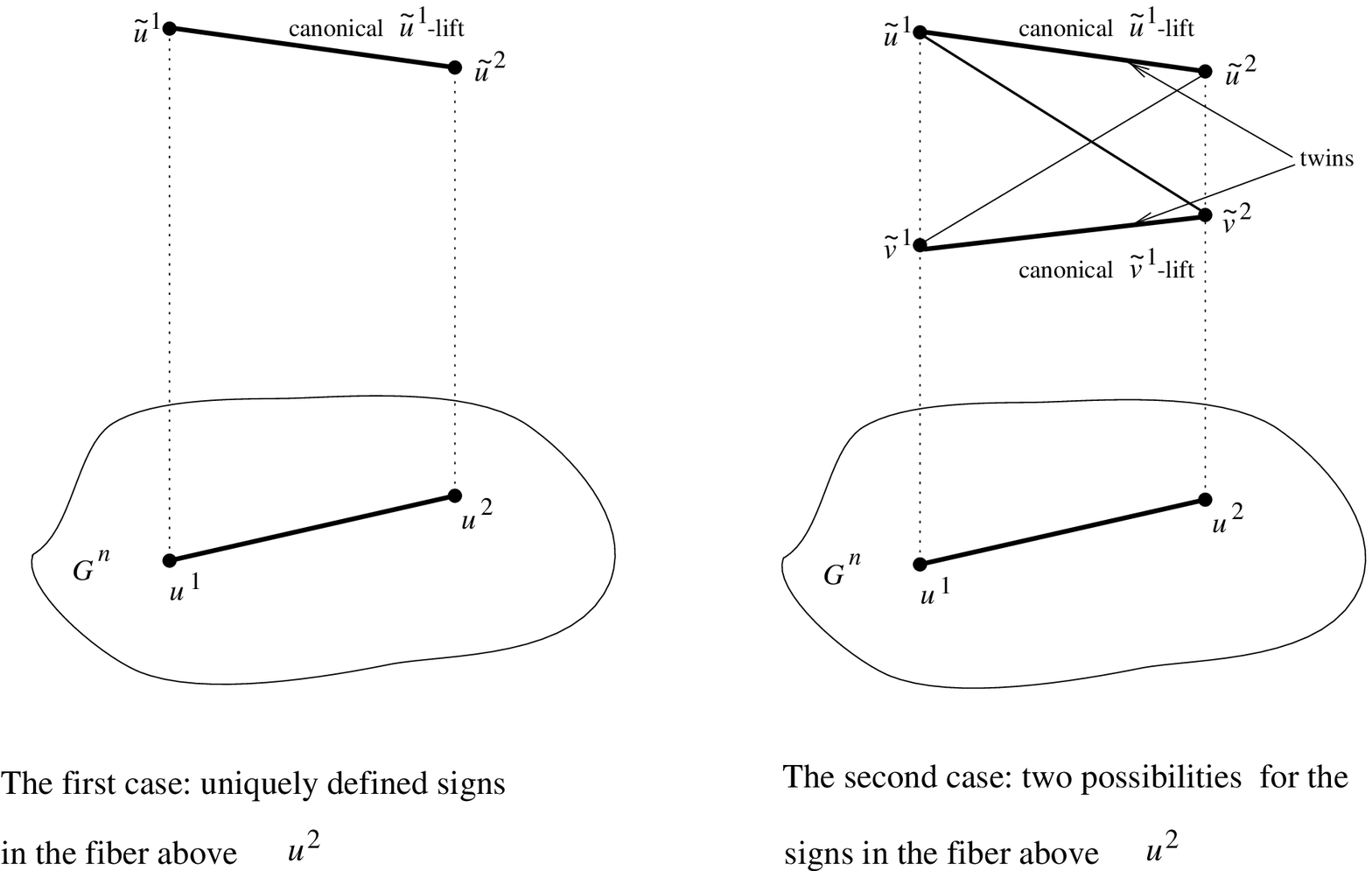}}}
\midspace{1mm}
\caption{Fig.5. Canonical lifts and twins}
\vskip 5pt

In any case, if $[\tilde u^1,\tilde u^2]$ is a canonical lift of $e$, we
write $\tilde u_2=c_e(\tilde u_1)$ and $\kap_2=c_e(\kap_1)$. We deliberately
drop the subscript if there is no ambiguity in identifying the edge in 
question.

Let $\gamma$ be a path in $G^n$. A path $\tilde \gamma$ in  $\widetilde G^n$
is said to be a canonical lift of $\gamma$ if all the edges of $\tilde\gamma$
are canonical lifts of edges of $\gamma$. Evidently, each path $\gamma$ 
in $G^n$ has a canonical lift, which depends on the choice of the initial 
vertex in $\widetilde G^n$. If we fix such a vertex, then the canonical lift
of $\gamma$ is defined uniquely.

\subheading{3.6. Special lifts} 
Consider a path $\gamma$ in $G^n$ such that some canonical
lift $\tilde e=[\tilde u,\tilde v]$ of its last edge has the twin 
$\tilde e'=[\tilde u',\tilde v']$. Let $\tilde \gamma$ and $\tilde \gamma'$
be the canonical lifts of $\gamma$ containing $\tilde e$ and $\tilde e'$,
respectively. We replace the edge $\tilde e'$ by the edge 
$[\tilde u,\tilde v']$ and get a noncanonical lift $\hat\gamma^2$ of the closed
path $\gamma^2$ in $G^n$ obtained by traversing $\gamma$ twice in the
opposite directions. The path $\hat \gamma^2$ is said to be a {\it special
lift\/} of $\gamma^2$ (see Fig.~6). Observe that the number of edges in a 
special lift is always even.

\vskip 10pt
\centerline{\hbox{\hskip 0.2cm \epsfxsize=15truecm\epsfbox{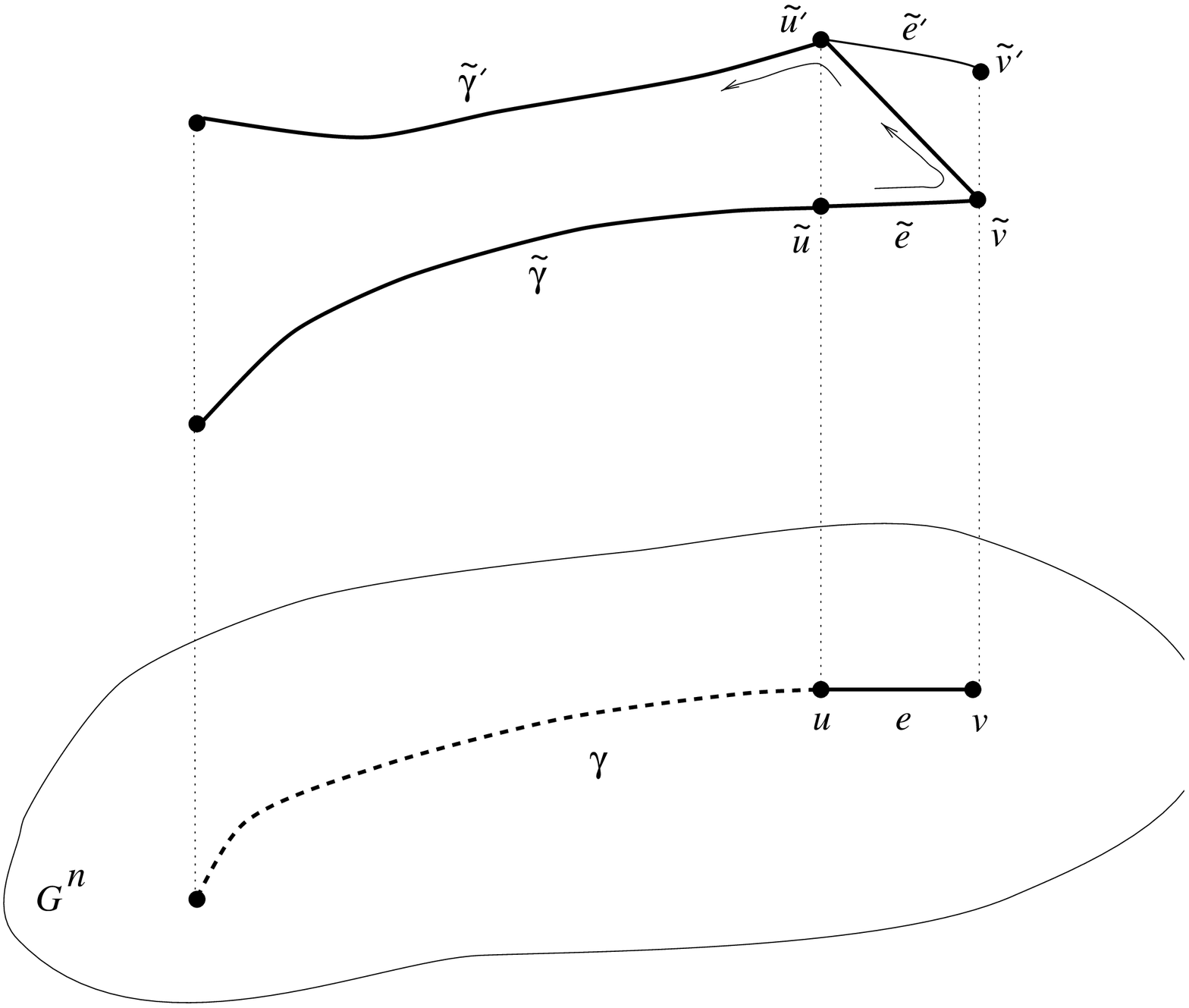}}}
\midspace{1mm}
\caption{Fig.6. The special lift of $\gamma^2$}
\vskip 5pt

\proclaim{3.7. Lemma} For any lift $\bar\gamma$ of a closed path $\gamma$ in
$G^n$ there exists a sequence $\bar\gamma_1,\dots,\bar\gamma_k$
of paths in $\widetilde G^n$ such that $\bar\gamma_k$ is a canonical lift 
of $\gamma${\rm,} all the other paths in the sequence are special lifts{\rm,}
and the concatenation of $\bar\gamma_1,\dots,\bar\gamma_k$ is a path between
the endpoints of $\bar\gamma$.
\endproclaim

\demo{Proof} Indeed, let us use the following procedure. Go along
$\bar\gamma$ until the first noncanonical edge occurs, then return
back using canonical edges only; we thus get a special lift.
Then traverse the second half of the previous path in the opposite direction
and continue along $\bar\gamma$ until the next noncanonical edge occurs. 
Return back using canonical edges only, and so on. The very last path in 
this sequence is a the canonical lift of $\gamma$.
\qed \enddemo

\subheading{3.8. Graphs $\Gamma^n(u)$ } Instead of $v_0^n$, one can start 
the construction procedure at any vertex $u$ of $G^n$, i.e., at any other 
arrangement of pseudolines. In this case the vertex set of $\Gamma^n(u)$ is 
defined as the preimage $\pi^{-1}(u)$. Obviously, the  pseudolines divide
the plane into a disjoint union of elementary regions (triangles and
diamonds in the case of $v_0^n$). Each elementary region has exactly two 
horizontal nodes, and one can define edges of $\Gamma^n(u)$ by exactly 
the same rule as above, with the help of admissible involutions.

Let $u$ and $v$ be two adjacent vertices of $G^n$. The edge $[u,v]$
corresponds to a 3-move that produces the 
``perestroika'' of the pseudoline arrangment $u$  defined by some elementary
triangle $A$, see Fig.~7.

\vskip 10pt
\centerline{\hbox{\hskip 0.2cm \epsfysize=1.5cm\epsfbox{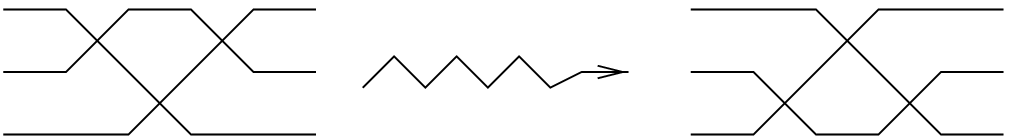}}}
\midspace{1mm}
\caption{Fig.7. Perestroika of diagrams corresponding to a 3-move} 
\vskip 5pt

We say that $A$ is the $e$-{\it core\/} of $u$.
Observe that there exists a  natural one-to-one correspondence between the
elementary regions of $u$ and $v$, which takes the $e$-core of $u$ to that 
of $v$; we denote this correspondence by $\hat c_e$ (as above, we drop the
subscipt when $e$ is identified unambigously).

\proclaim{3.9. Lemma} Let $e=[u,v]$ be an edge of $G^n${\rm,} $(u,\kap)\in
\pi^{-1}(u)$ be a vertex of $\widetilde G^n${\rm,} $A$ be the $e$-core
of $u${\rm,} and $B$ be an elementary region of $u$ such that $I_B$ is
admissible for $\kap$. Then exactly one
of the following possibilities holds\/{\rm:}

either $I_{\hat c(B)}$ is admissible for $c(\kap)$ and $c(I_B(\kap))=
I_{\hat c(B)}(c(\kap))${\rm;}

or $I_{\hat c(B)}$ is admissible for $c(\kap)${\rm,}
$I_{\hat c(A)}$ is admissible for $I_{\hat c(B)}(c(\kap))${\rm,}
and $c(I_B(\kap))=I_{\hat c(A)}\circ I_{\hat c(B)}(c(\kap))${\rm;}

else $I_{\hat c(A)}$ is admissible for $c(\kap)${\rm,}
$I_{\hat c(B)}$ is admissible for $I_{\hat c(A)}(c(\kap))${\rm,}
and $c(I_B(\kap))=I_{\hat c(B)}\circ I_{\hat c(A)}(c(\kap))$.
\endproclaim

\demo{Proof}
If $A$ and $B$ are not neighbors, then the first of the above possibilities
is apparently true.

All cases of neighboring $A$ and $B$ (up to obvious symmetries)
are considered in the following series of pictures. Abusing notation we 
omit $c$ and $\hat c$ to make figures more clear.

\newpage

1. Equilateral triangle (all the nodes of $A$ have the same sign).

Below we present all possible basic cases.

\vskip 10pt
\centerline{\hbox{\hskip 0.2cm \epsfysize=15cm\epsfbox{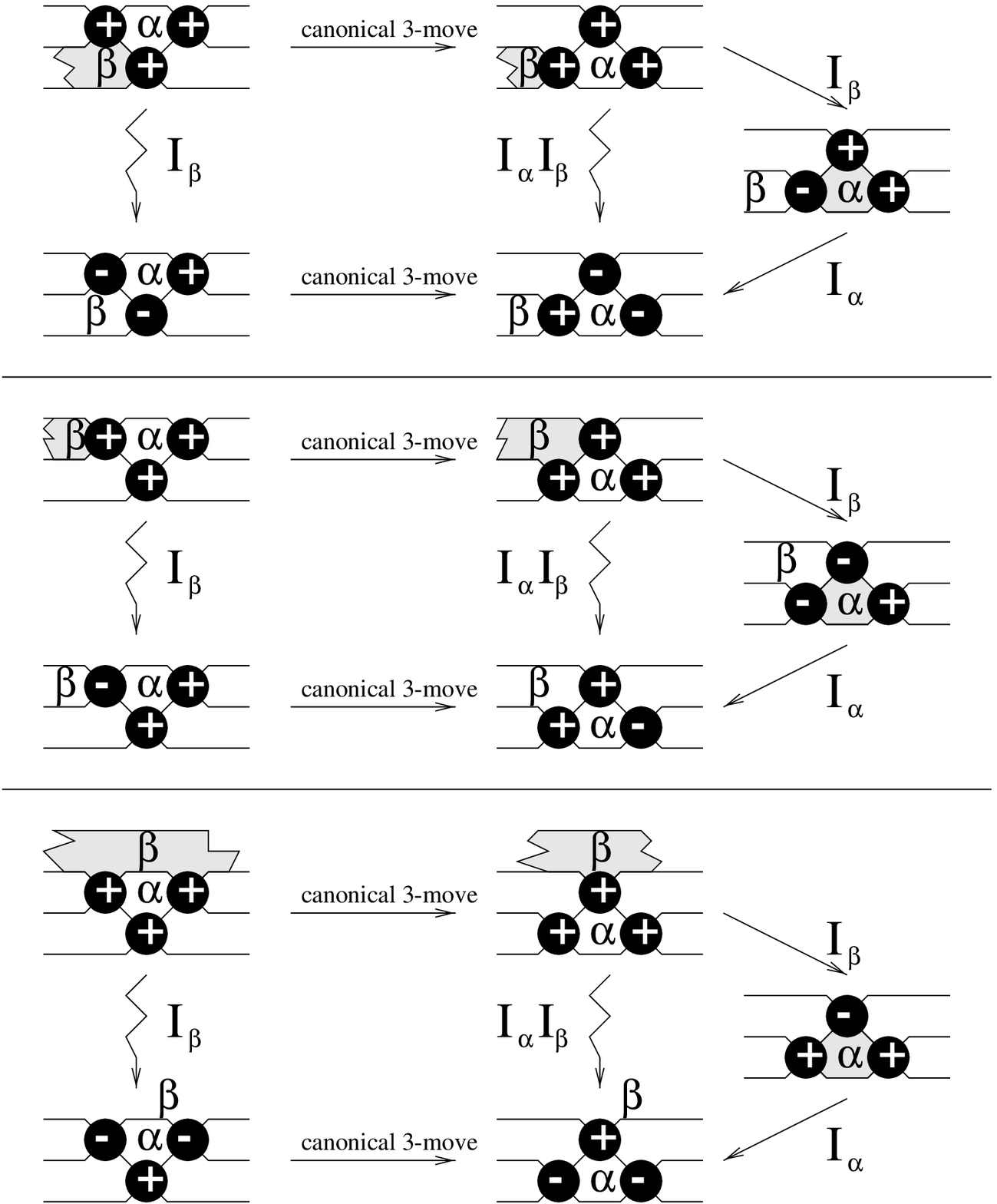}}}
\midspace{1mm}
\caption{Fig.8. Equilateral case}
\vskip 5pt

All other possible choices of signs and positions of $B$ are symmetric
(obtained by the global sign change and/or symmetry w.r.t. the vertical axis) 
to the ones considered.

\newpage

2. Isosceles triangle (base nodes of $A$ have the same sign).

Below we present all possible basic cases.

\vskip 10pt
\centerline{\hbox{\hskip 0.2cm \epsfysize=15cm\epsfbox{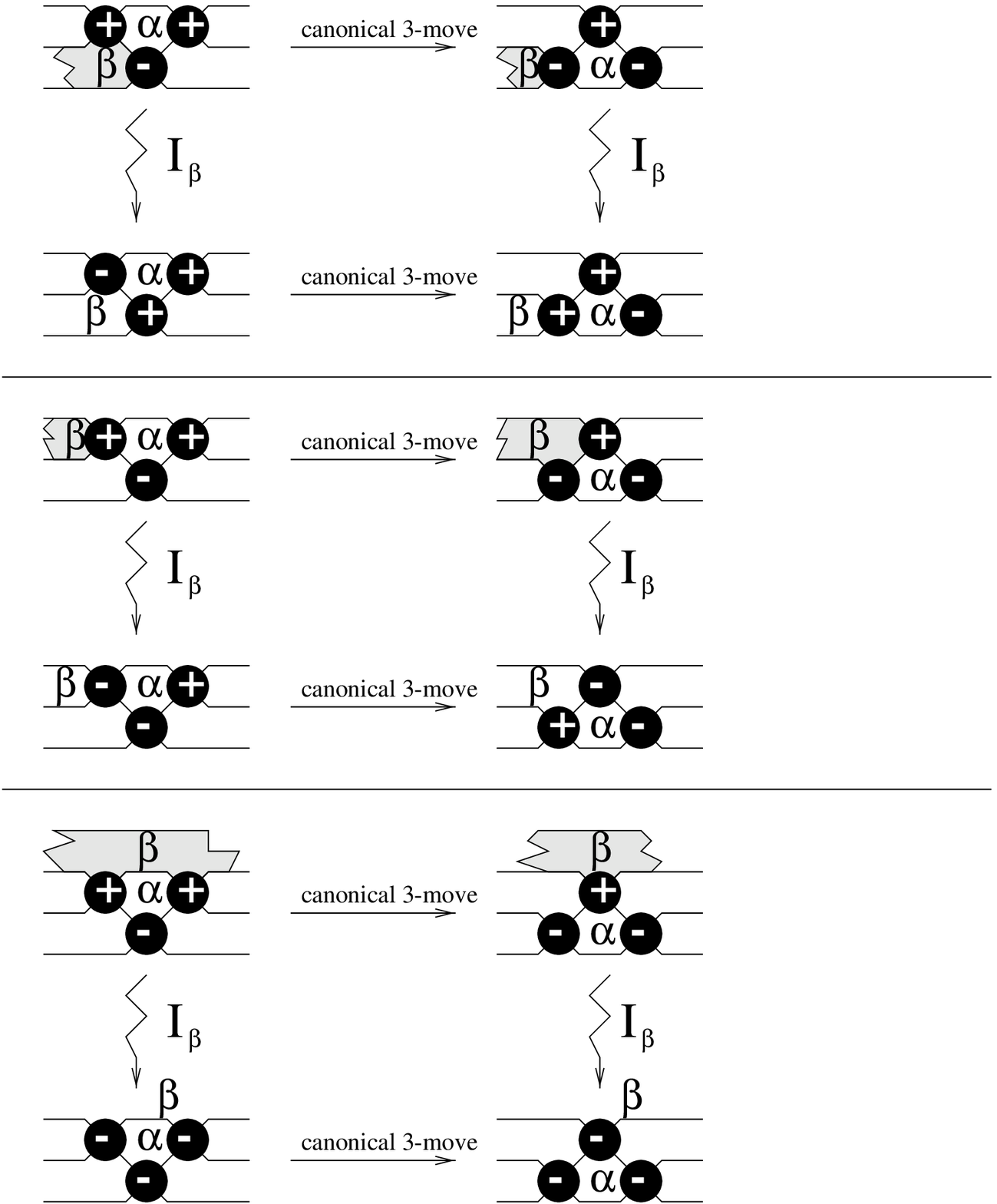}}}
\midspace{1mm}
\caption{Fig.9. Isosceles case}
\vskip 5pt

All other possible choices of signs and positions of $B$ are symmetric
to the ones considered.

\newpage

3. Scalene triangle (base nodes of $A$ have different signs).

Below we present all possible basic  cases.

\vskip 10pt
\centerline{\hbox{\hskip 0.2cm \epsfysize=16cm\epsfbox{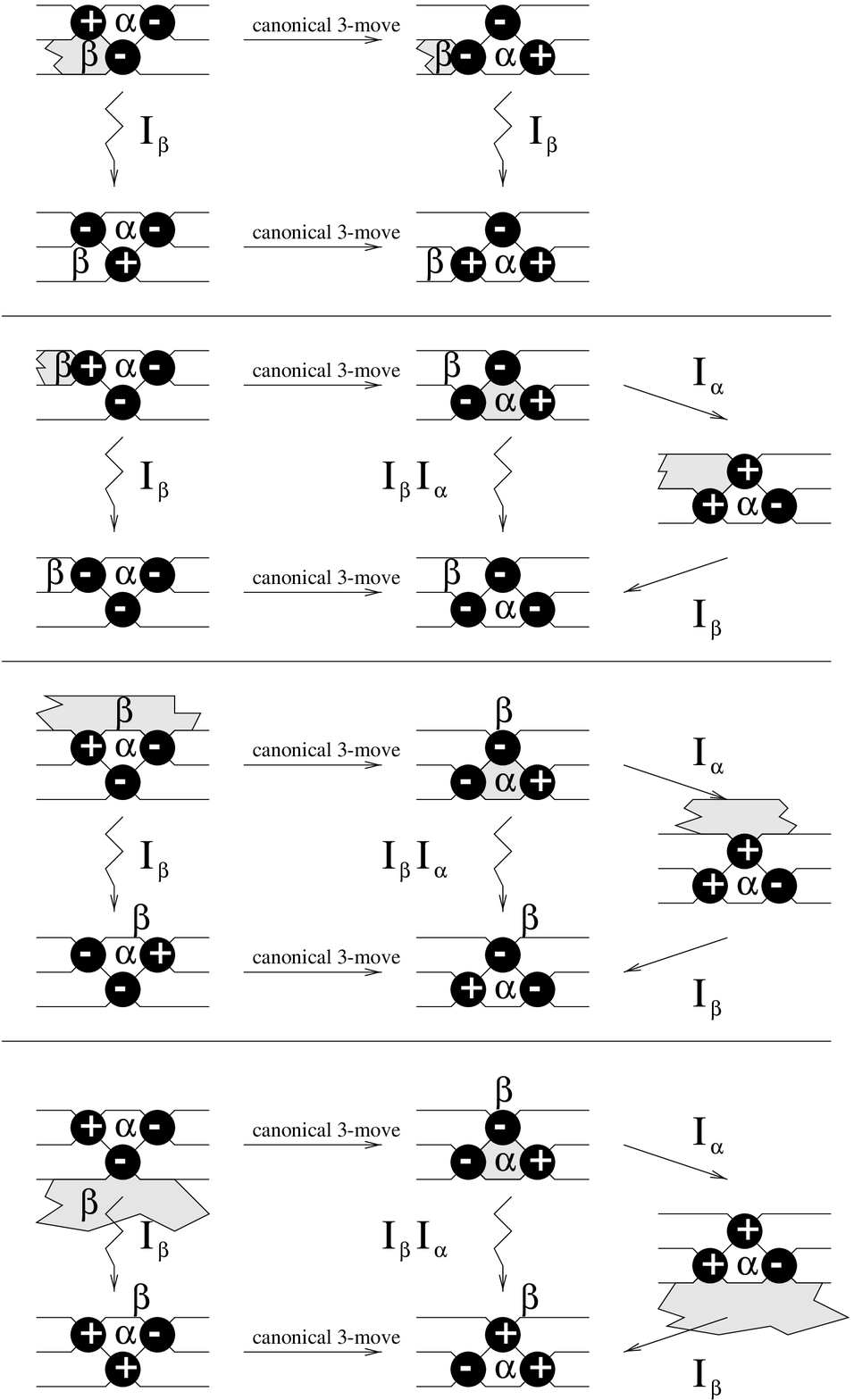}}}
\midspace{1mm}
\caption{Fig.10. Scalene case}
\vskip 5pt

All other possible choices of signs and positions of $\beta$ are symmetric to
the ones considered.

These series of pictures prove the lemma.
\qed
\enddemo

\proclaim{3.10. Corollary} The endpoints of a special lift lying in a fiber
$\pi^{-1}(u)$ for some $u\in G^n$ belong to the same 
connected component of $\Gamma^n(u)$. \endproclaim 

\demo{Proof}
We prove the statement by induction on $l$, where $l$ is the half-length
 of the special lift $\tilde\gamma$ in question.

 Let $l=1$, then $\tilde\gamma$ is a special lift of some 
edge $e$ and coincides (up to obvious symmetries) with the
path $\{a,b\}$ shown on Fig.~11. The existence of the $\Gamma^n(u)$-edge 
$c$ shown by the curved arrow proves the base of induction.

\vskip 10pt
\centerline{\hbox{\hskip 0.2cm \epsfxsize=10truecm\epsfbox{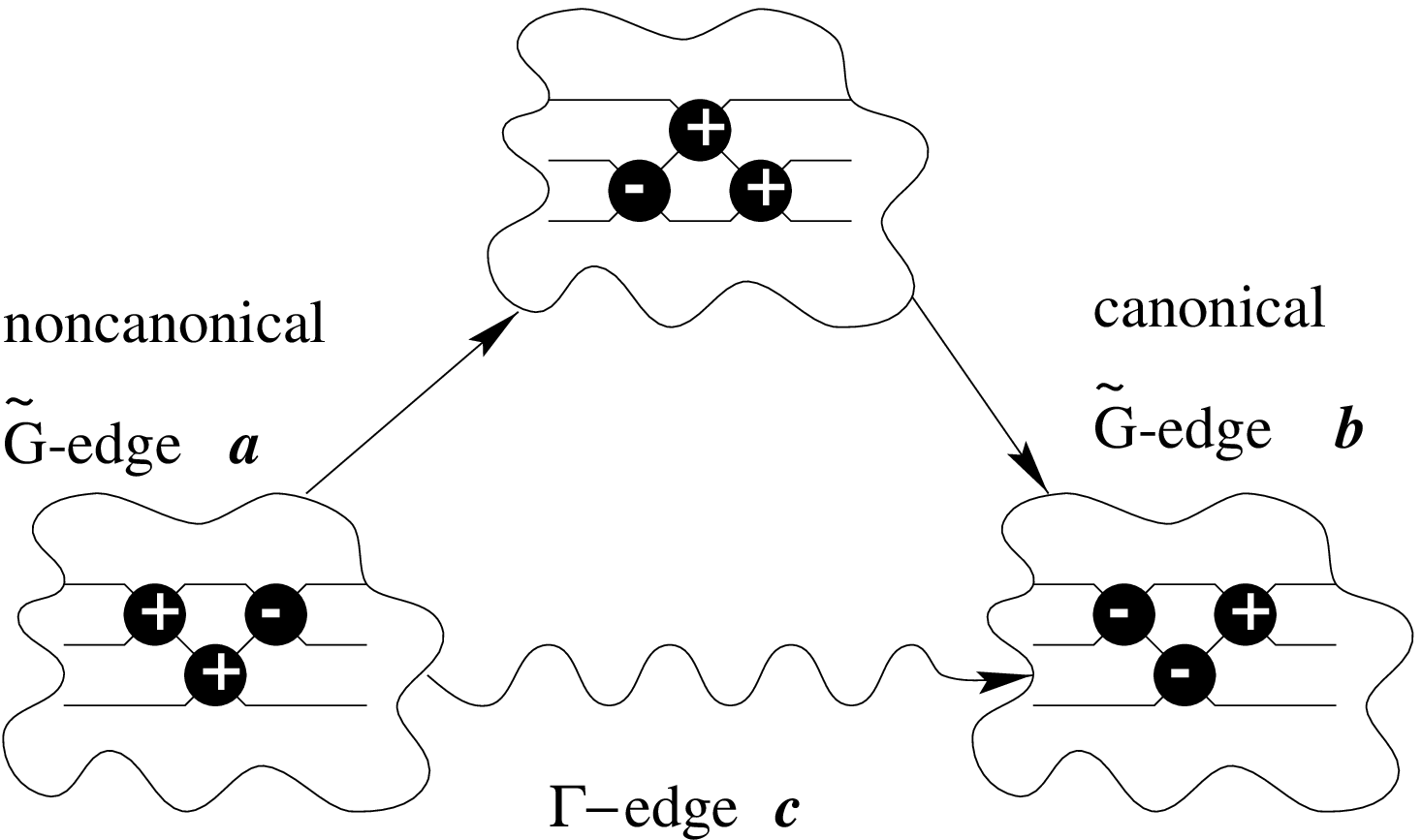}}}
\midspace{1mm}
\caption{Fig.11. Example of a special lift of length 2} \vskip 5pt

Suppose now that the statement holds for $l\ls K$. Let us
prove it for $l=K+1$. Consider a path 
$\gamma=(u=u_0,u_1,\dots,u_K,u_{K+1},u_K,\dots,u_1,u_0)$ in $G^n$ of
length $2 (K+1)$ and take its arbitrary
special lift $\tilde \gamma$ in $\tilde G$. Let $\tilde u_0$ and $\tilde u'_0$
be the endpoints of $\tilde\gamma$, $\tilde u_1=c_e(\tilde u_0)$, 
$\tilde u'_1=c_e(\tilde u'_0)$, where $e=[u_0,u_1]$. Then $\tilde u_1$ and
$\tilde u'_1$ are the endpoints of a special lift of length $2K$, and thus
they belong to the same connected component of $\Gamma^n(u_1)$. Now, 
by Lemma 3.9, $\tilde u_0$ and $\tilde u'_0$ belong to the same connected
component of $\Gamma^n(u)$. 
\qed\enddemo

\subheading{3.11} The following two types of cycles in $G^n$ are called
the 4-{\it cycle\/} and the 8-{\it cycle\/}, respectively, see Figs.~12,13. 

\vskip 10pt
\centerline{\hbox{\hskip 0.2cm \epsfxsize=10truecm\epsfbox{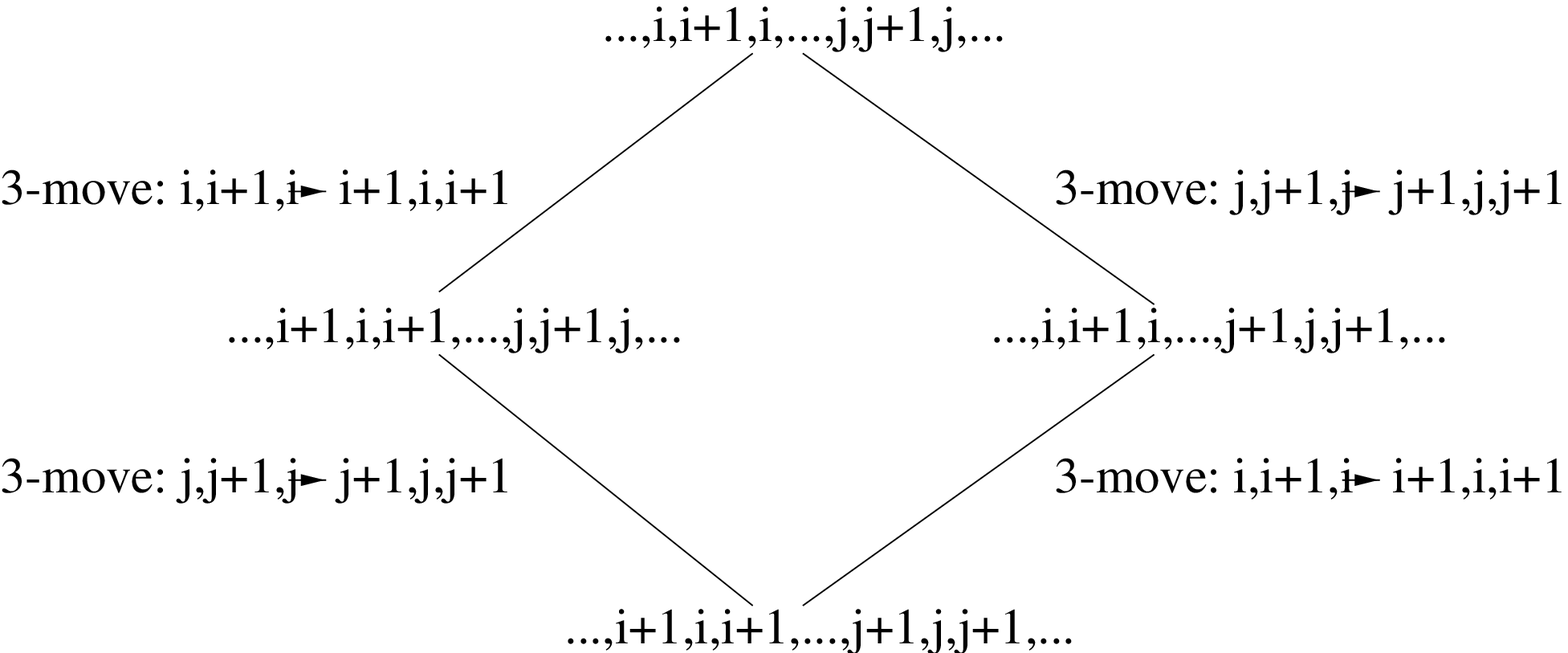}}}
\midspace{1mm}
\caption{Fig.12. A 4-cycle (commutativity of separated transpositions)} 
\vskip 5pt

\vskip 10pt
\centerline{\hbox{\hskip 0.2cm \epsfxsize=10truecm\epsfbox{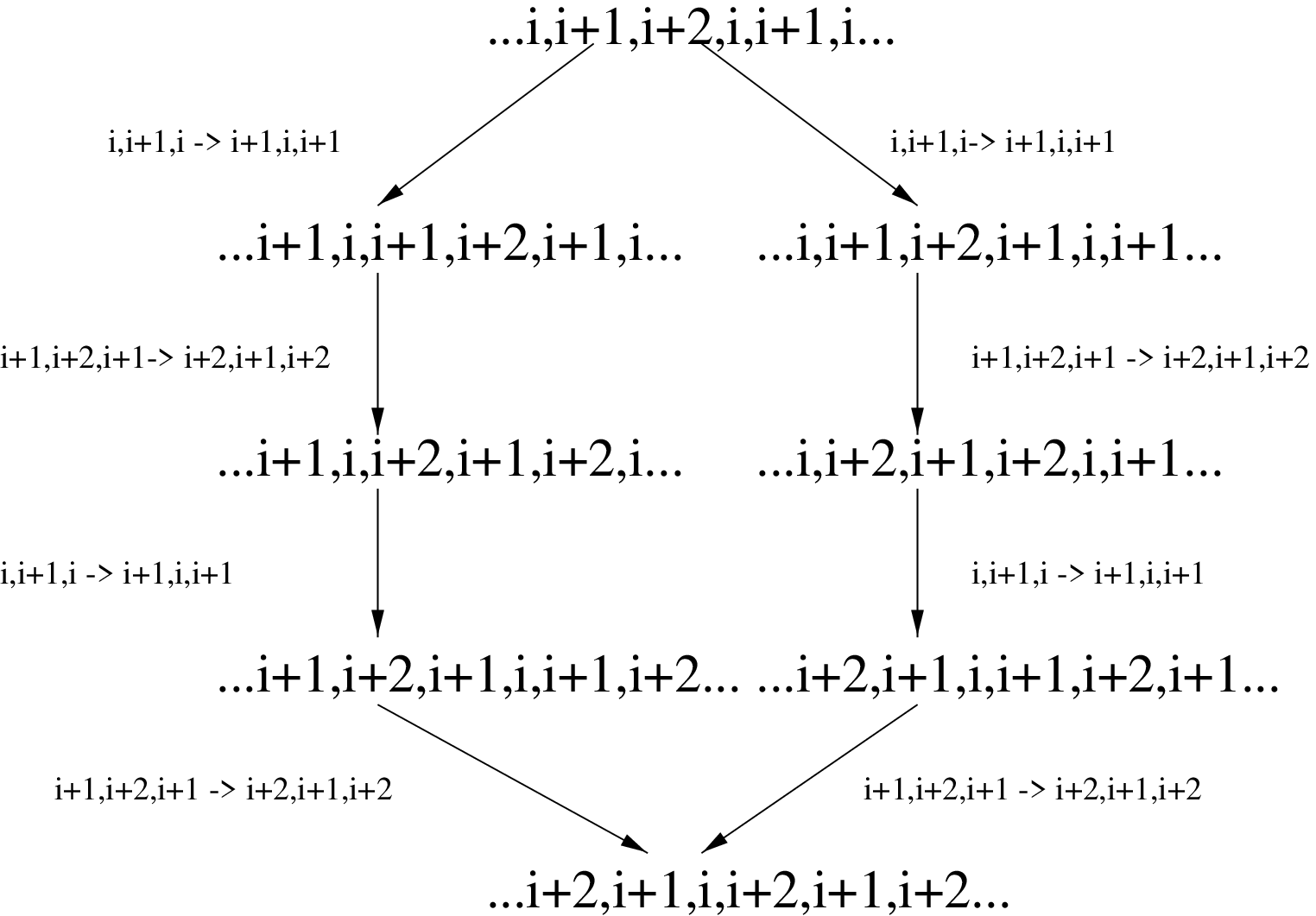}}}
\midspace{1mm}
\caption{Fig.13. An 8-cycle (Frenkel-Moore equation)} \vskip 5pt

\proclaim{3.12. Lemma} Any canonical lift of a $4$-cycle in $G^n$ is a cycle in
$\widetilde G^n$.
 \endproclaim

\demo{Proof} Follows from the commutativity property for 4-cycles
in $G^n$: the resulting sign transition for two consecutive nonintersecting
triples $\kap_i,\kap_{i+1},\kap_{i+2}$ and $\kap_k,\kap_{k+1},\kap_{k+2}$
does not depend on their order
\qed
\enddemo

For 8-cycles, the situation is  more complicated: a canonical lift of an
8-cycle in $G^n$ is not necessarily a cycle in $\widetilde G^n$. 
However, the following statement holds.

\proclaim{3.13. Lemma} The endpoints of a canonical lift of an $8$-cycle 
{\rm(}considered as a loop at $u\in G^n${\rm)}
belong to the same connected component of $\Gamma^n(u)$.
\endproclaim

\demo{Proof} One can prove this fact by a direct check of all the possibilities
occuring for diagrams with 4 pseudolines. The authors have also written a 
program in Mathematica that calculates the signs in the final diagram for all
possible sets of signs in the initial digram and checks if these vertices are 
connected in $\Gamma^n(u)$. 
\qed
\enddemo

\proclaim{3.14. Lemma} The set of all $4$- and $8$-cycles in $G^n$ form a 
system of generators for the first homology group $H_1(G^n,\Bbb Z_2)$.
\endproclaim

\demo{Proof} Let us introduce the rank of a reduced decomposition as 
follows: $\rho(s_{i_1}s_{i_2}\dots s_{i_k})=i_1+i_2+\cdots+i_k$.
Obviously, the rank does not change under 2-moves, thus it is well-defined 
for vertices of $G^n$. Any 3-move changes the rank of a vertex exactly by 1.

Let $\gamma$ be an arbitrary cycle in $G^n$, and denote by $\rho(\gamma)$ the
maximum of the ranks of the vertices in $\gamma$. We prove by induction 
on $\rho(\gamma)$ that $\gamma$ is homologous to a sum of 4- and 8-cycles. The
base of induction is trivial. Assume now that the statement holds for the
cycles with the rank less than $r$, and fails for cycles of rank $r$.
Let $\gamma$ be a counterexample  with the minimal number of vertices of
rank $r$. Consider an arbitrary vertex $v\in\gamma$ such that $\rho(\gamma)=r$,
and its neighbours $v'$ and $v''$ in $\gamma$. The edges $[v,v']$ and $[v,v'']$
correspond to 3-moves, each involving a triple of simple transpositions.
We distinguish the following three cases.

1. These triples coincide. Then $v'=v''$, and thus $\gamma$ is homologous
to a cycle $\gamma^*$ such that either $\rho(\gamma^*)=r-1$, 
or $\rho(\gamma^*)=r$ and the number of vertices of rank $r$ in $\gamma^*$
is less than that in $\gamma$, a contradiction.

2. These triples are disjoint. Then there exists a vertex $v^*$ of $G^n$
such that $(v^*, v',v, v'',v^*)$ from a 4-cycle $\gamma_4$; observe that
$\rho(v^*)=r-2$. Let $\gamma^*=\gamma+\gamma_4$, then either 
$\rho(\gamma^*)=r-1$, or $\rho(\gamma^*)=r$ and the number of vertices 
of rank $r$ in $\gamma^*$
is less than that in $\gamma$. Since $\gamma=\gamma^*+\gamma_4$, we thus 
see that $\gamma$ is not a counterexample.

3. These triples intersect, but not coincide. Then there exists an 8-cycle
$\gamma_8$ containing the edges $[v,v']$ and $[v,v'']$ such that the ranks
of all its vertices excluding $v$ are less than $r$. We then take $\gamma^*=
\gamma+\gamma_8$ and proceed as in the previous case.
\qed
\enddemo

\proclaim{3.15. Corollary} The endpoints of a canonical lift of any
 closed path in $G^n$ {\rm(}considered as a loop at $u\in G^n${\rm)}
belong to the same connected component of $\Gamma^n(u)$.
\endproclaim

\demo{Proof} Follows from Lemmas 3.9 and 3.12-3.14.
\qed
\enddemo

\proclaim{3.16. Proposition} If the endpoints of a path in $\widetilde G^n$ 
belong to $\Gamma^n${\rm,} then they lie in the same connected
component of $\Gamma^n$. \endproclaim

\demo{Proof} Follows from Lemma 3.7 and Corollaries 3.10 and 3.15.
\qed
\enddemo

\subheading{3.17}
{\it Proof of Theorem\/} 3.2. Follows directly from Propositions 3.3
and 3.16.
\qed 

\subheading{3.18. $\G_{n-1}$-action}
In order to interpret the obtained results as the $\bZ_2$-linear
action described in the introduction, we proceed as follows.
The nodes of $v_0^n$ form the uppertriangular shape (rotate Fig.~1 by $\frac
\pi 2$ ). Each vertex of $\Gamma^n$ corresponds uniquely to an
uppertriangular $(n-1)\times (n-1)$-matrix with $\bZ_2$-values (by convention,
$+$ corresponds to $1$ and $-$ to $0$.)
The action of the operators $I_A$ in terms of $\bZ_2$-matrices is
exactly the action of $\G_{n-1}$ presented in the introduction. Namely, 
if a matrix in $N^{n-1} (\bZ_2)$ has a dense $2\times 2$-submatrix
of one the following two forms:
$$
\pmatrix 1&*\\*&0\endpmatrix\text { or }\pmatrix 0&*\\*&1\endpmatrix,
$$ 
then the corresponding operator $I_A$ changes each entry in this submatrix
to the opposite one. If a submatrix is of the form
$$
\pmatrix 0&*\\*&0\endpmatrix\text { or }\pmatrix 1&*\\*&1\endpmatrix,
$$ 
then $I_A$ acts identically. In both cases the action of $I_A$ is just an 
addition of the trace of the submatrix to each entry. This gives exactly 
the action of $\G_{n-1}$.

\heading \S 4. Final remarks. \endheading
\subheading{4.1}
In this paper we have reduced the question on the number of connected 
components in the intersection of two open opposite Schubert cell in 
$SL_n(\bR)/B$ to 
the calculation of  the number of orbits  of  $\G_{n-1}$-action on the 
$\bF_2$-linear space $N^{n-1}(\bF_2)$,  see Introduction. 

Let us list some relatively obvious properties of $\G_n$.

a) All generators $g_{ij}$ of $\G_n$ are involutions and belong to the same 
conjugacy class.

b) Let us construct the following planar graph with the vertex set
$\{g_{ij}\}$, $1\ls i\ls j\le n-1$. We connect by edges all pairs 
of the form $(g_{ij},g_{i+1,j})$, $(g_{ij},g_{i,j+1})$ and
$(g_{ij},g_{i+1.j+1})$, i.e., we arrange the generators $g_{ij}$ into an 
equilateral triangle and place them on the hexagonal lattice in $\bR^2$ 
with edges joining the neighboring vertices of the lattice. Then

\qquad(i) any distant (i.e., not joined by an edge) generators commute;

\qquad(ii) any 2 generators $a$ and $b$ joined by an edge satisfy $(ab)^3=1$;

\qquad(iii) any 3 generators $a$, $b$, $c$ pairwise joined by edges generate
the group $S_4$.

c) The group $\G_n$ has a natural $\bF_2$-linear representation on
the space $N^{n-1}(\bF_2)$. Indeed, consider the linear map
$\pi\: N^n(\bF_2)\to N^{n-1}(\bF_2)$ such that the $(i,j)$th entry in the image
equals the sum of the $(i,j)$th and $(i+1,j+1)$th entries in the inverse image.
Since $\ker\pi$ is invariant under the action of $\G_n$ on $N^{n}(\bF_2)$,
one obtains the induced action of $\G_n$ on $N^{n-1}(\bF_2)$. It is easy to 
see that $g_{ij}$ acts by adding the $(i,j)$th entry of the matrix 
to all its neighbors in the hexagonal lattice.

d) There exists a scew-symmetric bilinear form $\Phi$ on $N^{n-1}(\bF_2)$
of corank $\left[\frac{n-1}{2}\right]$ that is preserved under the 
conjugate to the above $\G_n$-action on $N^{n-1}(\bF_2)$.

Note that a similar situation was  already studied earlier
in connection with the monodromy theory for isolated singularities,
see \cite{J1},\cite{J2},\cite{C},\cite{W}.

Although we do not have a complete group-theoretical description of $\G_n$,
we propose the following conjecture about its orbits in $N^n(\bF_2)$.
 
\proclaim{Conjecture} The $\G_n$-action on $N^n(\bZ_2)${\rm,} $n\gs5${\rm,}
have the following orbits\/{\rm:}

\noindent for $n=2k${\rm:} 

$2^{2k}$ orbits of length $1${\rm,}

$2^{k+2}(2^{k-1}-1)$ orbits of length $2^{(2k+1)(k-1)}${\rm,}

$2^k$ orbits of length $2^{(k+1)(k-1)}(2^{k(k-1)}-1)${\rm,}

$2^k$ orbits of length $2^{(k+1)(k-1)}(2^{k(k-1)}+1)${\rm,}

$2^{k+1}$ orbits of length $2^{k-1}(2^{2k(k-1)}-1)${\rm;}

\noindent for $n=2k+1${\rm:} 

$2^{2k+1}$ orbits of length $1${\rm,}

$2^{k+2}(2^{k}-1)$ orbits of length $2^{(2k-1)(k+1)}${\rm,}

$2^{k+1}$ orbits of length $2^{k^2+k-1}(2^{k^2}-1)${\rm,}

$2^{k+1}$ orbits of length $2^k(2^{k^2}-1)(2^{k^2-1}+1)$.
\endproclaim

As an immediate corollary we get the proof of the main conjecture stated 
in the introduction.

\subheading{4.2. The case of nonopposite flags $f$ and $g$}
It is well known that the orbits of the $SL_n$-action on the space of 
pairs of flags are parametrized by permutations. One thus can ask for
finding the number of connected components in the intersection $U^n_{f,g}$ 
of open Schubert cells w.r.t. flags $f$ and $g$ provided
the flags $f$ and $g$ are in relative position $w\in S_n$. In the particular
case $w=(n\;n-1\;\dots 1)$ we get our initial problem concerning opposite
cells. 

We define an {\it affine pseudoline arrangement} similar to usual pseudoline
arrangement with the only distinction: we allow parallel pseudolines (that is,
intersecting at infinity).  Compact connected components of the complement of 
an affine pseudoline arrangement we call {\it elementary regions}.

Given a permutation $w$, let us consider any 
affine pseudoline arrangement $\P$ realizing $w$ in the obvious way, 
cf. \cite{BFZ}. Let $R$ denote the set of intersection points in $\P$, 
$W$ be the $\bF_2$-vector space with the basis $R$, and $V=W^*$. 
For each elementary region $a$ we define a linear operator $g_a\: V\to V$ 
in the following way. Each elementary region has exactly 2 ``horizontal'' 
boundary vertices. Let $p_1,p_2,\dots, p_q$ be boundary vertices of $a$ and 
$p_1$, $p_2$ be ``horizontal''. For any intersection point $p\in R$ we 
denote by  $\psi_p\in V$ the characteristic function of the point $p$.
Put
$$
g_a(f)=f+(f(p_1)+f(p_2))(\psi_{p_1}+\psi_{p_2}+\dots+\psi_{p_q}).
$$

Denote by $\G(w,\P)$ the subgroup of $GL(V)$ generated by all $g_a$.
The following statement would be a generalization of the main result of this
paper.

\proclaim{Conjecture} The number of connected components in $U^n_{f,g}$
is equal to the number of orbits of $\G(w,\P)$-action in $V$.
\endproclaim

\enddocument

In this paper we reduce the question on the number of connected components 
in the intersection of 2 open opposite Schubert cell in $SL_n/B$ to 
the enumeration of equivalence classes of a purely
combinatorial flip action on the set $N^{n-1}(\bZ_2)$ of upper 
triangular $(n-1)\times (n-1)$-matrices with $\bZ_2$-valued entries, 
see Introduction. The restriction that
 the sign reversion is only permitted if the diagonal entries have opposite
signs makes the problem of finding the number of equivalence classes 
complicated (and interesting).
The authors have an explicit conjecture about the number and the structure
of these equivalence classes, as well as proofs of several related 
propositions. 

\proclaim{Conjecture} The flip action on $N^n(\Bbb Z_2)${\rm,} $n\gs5${\rm,}
have the following equivalence classes\/{\rm:}

\noindent for $n=2k${\rm:} 

$2^{2k}$ classes of length $1${\rm,}

$2^{k+2}(2^{k-1}-1)$ classes of length $2^{(2k+1)(k-1)}${\rm,}

$2^k$ classes of length $2^{(k+1)(k-1)}(2^{k(k-1)}-1)${\rm,}

$2^k$ classes of length $2^{(k+1)(k-1)}(2^{k(k-1)}+1)${\rm,}

$2^{k+1}$ classes of length $2^{k-1}(2^{2k(k-1)}-1)${\rm;}

\noindent for $n=2k+1${\rm:} 

$2^{2k+1}$ classes of length $1${\rm,}

$2^{k+2}(2^{k}-1)$ classes of length $2^{(2k-1)(k+1)}${\rm,}

$2^{k+1}$ classes of length $2^{k^2+k-1}(2^{k^2}-1)${\rm,}

$2^{k+1}$ classes of length $2^k(2^{k^2}-1)(2^{k^2-1}+1)$.
\endproclaim

As an immediate corollary we get the proof of the main conjecture stated 
in the introduction.

\enddocument